\documentclass[reqno]{amsart}
\usepackage{amssymb,amsmath,amsthm,amstext,amsfonts}
\usepackage[dvips]{graphicx}
\usepackage{psfrag}
\usepackage{url}
\usepackage{amsmath,amstext,amsthm,amsfonts}
\usepackage{color}
\usepackage{xcolor}

\usepackage[colorlinks=true, linkcolor=blue, urlcolor=red, citecolor=blue, hyperindex, backref]{hyperref}
\usepackage{ mathrsfs }

\makeatletter \@addtoreset{equation}{section} \makeatother

\renewcommand\thetable{\thesection.\@arabic\c@table}

\theoremstyle{plain}
\newtheorem{maintheorem}{Theorem}

\newtheorem{maincorollary}{Corollary}

\newtheorem{mainquestion}{Question}

\newtheorem{theorem}{Theorem }[section]
\newtheorem{proposition}[theorem]{Proposition}
\newtheorem{lemma}[theorem]{Lemma}
\newtheorem{corollary}[theorem]{Corollary}

\newtheorem{remark}[theorem]{Remark}

\newtheorem{definition}[theorem]{Definition}
\newtheorem{problem}{Problem}

\newtheorem{question}{Question}

\newcommand{\Sc}{\mathbb{S}}
\newcommand{\Lo}{\mathcal{L}}

\newcommand{\C}{\mathbb{C}}

\newcommand{\N}{\mathbb{N}}
\newcommand{\R}{\mathbb{R}}
\newcommand{\T}{\mathbb{T}}

\newcommand{\cL}{\mathcal{L}}

\def\ds{\displaystyle}

\begin{document}

\title{Density of spectral gap property for positively expansive dynamics and smooth potentials, with applications to the phase transition problem}

%\title{Thermodynamic formalism for positively expansive dynamics and smooth potentials}

\author{ Thiago Bomfim and  Victor Carneiro}

\address{Thiago Bomfim, Departamento de Matem\'atica, Universidade Federal da Bahia\\
Av. Ademar de Barros s/n, 40170-110 Salvador, Brazil.}
\email{tbnunes@ufba.br}
\urladdr{https://sites.google.com/site/homepageofthiagobomfim/}

\address{Victor Carneiro, Departamento de Ciências Extasas, Universidade Estadual de Feira de Santana\\
Av. Transnordestina, s/n - Novo Horizonte,
44036-900, Feira de Santana, Brazil.}
\email{victor.carneiro93@gmail.com}

\date{\today}

\begin{abstract}
It is known that all uniformly expanding dynamics $f: M \rightarrow M$ have no phase transition with
respect to a Hölder continuous potential $\phi : M \rightarrow \R$, in other words, the topological pressure function $\R \ni t \mapsto P_{top}(f , t\phi)$ is analytical. Moreover, the associated transfer operator $\mathcal{L}_{f , t\phi}$, acting on the space of Hölder continuous functions, has the spectral gap property for all $t \in \R$. For dynamics that are topologically conjugate to an expanding map, a full understanding has yet to be achieved. On the one hand, by \cite{KQW21,KQ22}, for such maps and continuous potentials, the associated topological pressure function can behave wildly. On the other hand, by \cite{BF23}, for transitive local diffeomorphisms on the circle and a large class of H\"older continuous potentials, the phase transition does not occur, and the associated transfer operator  has the spectral gap property for all parameters $t \in \R$. As a first approach to understanding what happens in high dimensions, we study positively expansive local diffeomorphisms in this paper. In particular, we show that the
associated transfer operator has the spectral gap property for a large class of regular potentials. Moreover, for a class of intermittent skew-products and a large class of regular potentials, we obtain phase transition results analogous to \cite{BF23}.
\end{abstract}

%\tableofcontents

\maketitle

%%%%%%%%%%%%%%%%%%%%%%%%%%%%%%%%
\section{Introduction}

From a dynamical systems point of view, phase transitions are often associated with the non-uniqueness of equilibrium states or lack of regularity of the topological pressure function, as a function of the potential.

More precisely, given $f : \Lambda \rightarrow \Lambda$ a continuous dynamical system on a compact metric space $\Lambda$ and $\phi : \Lambda \rightarrow \R$ a continuous potential, the variational principle for the pressure asserts that
$$
P_{top}(f, \phi) = \sup\{h_{\mu}(f) + \int\phi d\mu : \mu \text{ is an } f-\text{invariant probability}\}
$$
where $P_{top}(f, \phi)$ denotes the topological pressure of $f$ with respect to $\phi$ and $h_{\mu}(f)$ denotes the Kolmogorov-Sinai metric entropy. An equilibrium state $\mu_{\phi}$ for $f$ with respect to $\phi$ is a probability that attains the supremum. In the special case where $\phi \equiv 0$, the topological pressure $P_{top}(f, \phi)$ coincides with the topological entropy $h_{top}(f)$, which is one of the most important topological invariants in dynamical systems (see e.g. \cite{W82}).

In this paper, we say that $f$ has a (thermodynamic) phase transition with respect to $\phi$ if the topological pressure function
$$
\R \ni t \mapsto P_{top}(f, t\phi)
$$
is not analytic, we also say that $\phi$ has a phase transition.

The topological pressure function is expected to be regular. Indeed, it is known from Ergodic Theory that if $f$ has finite topological entropy, then the topological pressure function $\R \ni t \mapsto P_{top}(f, t\phi)$ is Lipschitz continuous. Furthermore, it follows from Rademacher's Theorem that $\R \ni t \mapsto P_{top}(f, t\phi)$ is differentiable almost everywhere. By \cite{W92}, if $f$ is expansive and has finite topological entropy, then the loss of differentiability of the topological pressure function is related to the non-uniqueness of the equilibrium states associated with the potential $\phi$. Therefore, the phase transition problem is related to such inconsistencies of thermodynamic quantities.

The phase transition problem is well understood in the classical case for expanding maps; however, only partial results have been obtained for other non-uniformly hyperbolic dynamics. By \cite{S72, Bow75, BR75}, all transitive hyperbolic or expanding dynamics $f$ have no phase transition with respect to any H\"older continuous potential $ \phi $. On the other hand, several non-uniformly hyperbolic examples are known to exhibit phase transitions with respect to regular or geometric potentials. For example:
\begin{itemize}
\item Manneville-Pomeau maps and geometric potential \cite{Lo93},
\item a large class of interval maps with an indifferent fixed point and geometric
potential \cite{PS92},
\item certain quadratic maps and geometric potential \cite{CRL13, CRL15, CRL19},
\item certain smooth maps of the interval and geometric potential \cite{CRL21},
\item porcupine horseshoes and geometric potential \cite{DGR14},
\item geodesic flow on  Riemannian non-compact manifolds with variable pinched negative sectional curvature and suitable H\"older continuous potential \cite{IRV18},
\item geodesic flow on certain  $M-$puncture spheres and geometric potential \cite{V17}.
\end{itemize}
Therefore, a complete characterization of the set of dynamics exhibiting phase transitions is an open problem.

Additionally, by \cite{KQW21, KQ22}, for some Markovian dynamics and
continuous potentials, the associated topological pressure function can behave
wildly. Thus, regularity of the potentials is required to make reasonable predictions. 

In view of this discussion, Bomfim and Carneiro \cite{BC21} proposed the following problem:

\begin{problem}\label{probA}
Which mechanisms are responsible for the occurrence of thermodynamic phase transitions for  $C^{1}$ local diffeomorphisms with positive topological entropy and H\"older continuous  potentials?
\end{problem}

Notably, significant progress was initiated in \cite{BC21,BF23}, answering this question for the case of circle maps.

When $f$ is a mixing expanding map and $\phi$ is a suitable potential, the analyticity of the pressure function is classically established via the Ruelle-Perron-Frobenius or transfer operator $\mathcal{L}_{\phi}$. This operator acts on functions $g : M \rightarrow \C$ and is defined as follows:
$$
\mathcal{L}_{f,\phi}(g)(x) := \sum_{f(y) = x}e^{\phi(y)}g(y). $$
More specifically, using the fact that $\mathcal{L}_{\phi}$ has the spectral gap property acting on a suitable Banach space, it is shown that the pressure function is analytic and, therefore, $f$ has no phase transition with respect to suitable potentials (see e.g. \cite{PU10}).

\begin{definition}
Let $E$ be a Banach space and $T : E \rightarrow E$ be a bounded linear operator. We say that $T$ has the {\bf spectral gap property} if there exists a decomposition of its spectrum $sp(T) \subset \C$ of the form
$sp(T) = \{\lambda_{1}\} \cup \Sigma_{1}$, where $\lambda_{1} >0$ is an isolated simple eigenvalue and there exists $0 < \lambda_{0} < \lambda_{1}$ such that $\Sigma_{1} \subset \{z \in \C : |z| < \lambda_{0}\}$. In this case, $\lambda_1 = \rho(T)$ is the spectral radius of $T$.
\end{definition}

Remarkably, the spectral gap property is extremely useful in obtaining fine statistical and differentiability properties of thermodynamic quantities, including equilibrium states, mixing behaviour, large deviations and topological pressure stability. (see e.g. \cite{Ba00,GL06, BCV16,BC19}).

 Employing the properties of transfer operators, Bomfim and Fernandes \cite{BF23} 
 showed that for transitive local diffeomorphisms of the circle, the absence of phase transitions (or, equivalently, the occurrence of the spectral gap property for the transfer operator acting on the space of H\"older continuous functions) is a typical property with respect to potentials. More formally:

\begin{theorem}\cite[Corollary A and Theorem B]{BF23}\label{comonBF}
Let $f : \Sc^{1} \rightarrow \Sc^{1}$ be a transitive and non invertible $C^{1}-$local diffeomorphism. Then:
\begin{enumerate}
    \item $\{ \phi : \Sc^{1} \rightarrow \R \text{ smooth } ;\; \phi \text{ has no  phase transition }\}$ is dense in the uniform topology;\\
    \item $\{\phi  : \Sc^{1} \rightarrow \R \text{ H\"older continuous 
    %and with bounded variation
    } ;\;
    \phi \text{ has phase transition }\}$ is not dense in the uniform topology.
\end{enumerate}
Moreover, given $\phi : \Sc^{1} \rightarrow \R$ an $\alpha-$H\"older continuous potential,
%and has bounded variation
then are equivalent:

\begin{enumerate}
    \item $\phi$ has no thermodynamic phase transition;\\
    \item $\mathcal{L}_{f, t\phi}$ has spectral gap property, acting on the $\alpha-$H\"older continuous functions space, for all $t \in \R$;\\
\end{enumerate}
\end{theorem}

This previous result answers the following question in the one-dimensional case:

\begin{question}\cite[Question E]{BC21}\label{quest34}
Let $f : M\rightarrow M$ be a local diffeomorphism $C^{1}$ on a manifold $M$ with $h_{top}(f) > 0$. The set of potentials $\phi$, in a suitable Banach space (smooth or H\"older continuous functions), such that $\mathcal{L}_{f, \phi}$ does not have the spectral gap property acting on a suitable Banach space (smooth or H\"older continuous functions) can be dense or residual in the uniform topology?
\end{question}

In this paper, we want to solve the previous problem for positively expansive dynamics. In particular, we will obtain a partial extension of Theorem \ref{comonBF}.

This paper is organized as follows. In Section \ref{sec:setmain}, we provide some definitions and a statement of the main results. In Section \ref{sec:preli} we recall the necessary framework on ergodic theory and the transfer operator.
In Section \ref{sec:proof}, the main results are proved. Finally, additional questions are proposed in Section \ref{sec:quest}.

\break 

\section{Definitions and statement of the main results}\label{sec:setmain}

This section is devoted to providing the setting and the statement of the main results.

\subsection{Setting}\label{subsec:set}

 Throughout the article, we will consider (positively) expansive $C^{r}-$local diffeomorphisms $f : M \rightarrow M$, where $M$ is a connected and compact Riemann manifold and $r \geq 1$ is an integer. Additionally, we ask that $f$ has {\it at least one repeller periodic point}.
 Remember:

\begin{definition}
    We say that $f$ is {\bf (positively) expansive} if there is a constant $\epsilon_{0} > 0$, called the expansivity constant, such that for every pair of points $x \neq y \in M$ there is $n \in \N$ such that $d(f^{n}(x),f^{n}(y)) \geq \epsilon_{0}.$
\end{definition}

\begin{definition}
    We say that $p \in M$ is a {\bf repeller periodic point} for $f$ if there is $k \in \N$ such that $f^{k}(p) = p$ and $||Df^{k}(p) \cdot v|| > ||v||$, for all $v \in T_{p}M\setminus \{0\}.$
\end{definition}

Observe that, in our setting, $f$ will be topologically conjugated to a uniformly expanding map (for more details, see \cite{CR80}).

%\begin{example}\label{examp}
 %Let $g: \T^d \to \T^d$ be expansive and $\alpha_1, \dots, \alpha_k \in \Sc^1$. Take $F: \T^d \times \Sc^1 \to \T^d \times \Sc^1$ a $C^r$-local diffeomorphism defined by $F(x,y) = (g(x), f_x(y))$, where $|f'_x(y)|>1$ for every $y\neq \alpha_j$, $f_{x|[\alpha_j,\alpha_{j+1})}$ is injective and  $f_x([\alpha_j,\alpha_{j+1})) = \Sc^1$ for all $j = 1, \ldots k-1$. So $F$ is expansive.
%\end{example}

%We will consider two particular cases:

%\begin{example}
%Let $F$ be a $C^{r}-$local diffeomorphism as in the previous example. We will say that $F \in TM1$ if $g$ is expanding. We will say that $F \in TM2$ if $g$ is intermittent, in other words, $d = 1$ and there exists  $x_1, \dots, x_l \in \Sc^{1}$ such that $|g'(x)|>1$ for every $x\neq x_j$, $g|_{[x_j,x_{j+1})}$ is injective and  $g([x_j,x_{j+1})) = \Sc^1$ for all $j = 1, \ldots l-1$.
%\end{example}

We want to construct interesting examples of local diffeomorphisms in higher dimensions that satisfy our hypotheses. Namely, we consider skew-products of expansive torus maps and intermittent circle maps. Let us define intermittent maps:

\begin{definition}
   We say that a map $f:\Sc^1 \to \Sc^1$ is {\bf almost expanding} if $f$ is transitive,  $|f'|\geq 1$ and there exist a finite number of fixed points $p_{1}, \dots, p_{k} \in \Sc^1$ such that if $|f'(x)|=1$ then $x = p_{j}$, for some $j = 1, \ldots, k$. We will say that $f$ is {\bf intermittent} if $f$ is almost expanding but $f$ is not expanding.
\end{definition}

Classic examples of intermittent maps are the Manneville-Pomeau maps. Our main examples will be in the following form:

\begin{definition}\label{examp}
  We say that a $C^{r}-$local diffeomorphism $F: \T^d \times \Sc^1 \to \T^d \times \Sc^1$ is an \textbf{almost expanding skew-product} if it is defined by $F(x,y) = (g(x), f_x(y))$, with $g: \T^d \to \T^d$ expansive and each $f_x:\Sc^1 \to \Sc^1$ is an almost expanding map with fixed points $q_{1,x}, \dots, q_{k,x} \in \Sc^1$ (dependent on $x$). 
\end{definition}

Note that an almost expanding skew-product is an expansive $C^{r}-$local diffeomorphism. We extend our main results to some toy models considered in the following.

\begin{definition}
Let $F:\T^d\times\Sc^1\to \T^d\times\Sc^1$ be an almost expanding skew-product. We will say that:
\begin{itemize}
\item $F \in TM1$ if $g$ is expanding and $\T^{d} \ni x \mapsto q_{l,x} \in \Sc^{1}$ is Lipschitz, for all $l = 1 ,\ldots, k$;
%$C^{1}$ and $\ds||Dg(x)|| > 1 + \frac{d q_{l,x}(x)}{d x},$ for all $x \in \T^{d}$ and $l = 1 ,\ldots, k$;
\item[]
\item $F \in TM2$ if $g$ is expanding  and the fixed points $q_{1,x}, \dots, q_{k,x}$ are constant with respect to $x$;
\item[]
\item $F \in TM3$ if $g$ is intermittent and the fixed points $q_{1,x}, \dots, q_{k,x}$ are constant with respect to  $x$. 
\end{itemize}
\end{definition}

\subsection{Main results}

The first result guarantees the density of the spectral gap property:

\begin{maintheorem}\label{theorA}
    Let $f : M \rightarrow M$ be an expansive $C^{r}-$local diffeomorphism on the connected and compact manifold $M$. Suppose that $f$ has at least one repeller periodic point. Then, there exists $\mathcal{H} \subset C(M , \R)$ open and dense subset such that: if $\phi \in \mathcal{H} \cap C^{r}(M , \R)$ then $\mathcal{L}_{f,\phi}$ has spectral gap property, acting on $C^{r}(M ,\C)$, and there exist $t_{0} , t_{1} > 0$ such that $\mathcal{L}_{f,t\phi}$ has spectral gap property, acting on $C^{r}(M ,\C)$, for all $|t| \geq t_{0}$ or $|t| \leq t_{1}$. 
\end{maintheorem}

In particular;

\begin{maincorollary}\label{maincorA}
Let 
$$X := \ds \Big\{\phi \in C^{r}(M , \R) : \mathcal{L}_{f,\phi|C^{r}} \text{ has spectral gap property }\Big\} \cap
$$
$$\Big\{\phi \in C^{r}(M , \R) : \exists t_{0} > 0 \text{ such that } \mathcal{L}_{f,t\phi|C^{r}} \text{ has spectral gap property } \forall |t| \geq t_{0} \Big\} \cap$$
$$\Big\{\phi \in C^{r}(M , \R) : \exists t_{1} > 0 \text{ such that } \mathcal{L}_{f,t\phi|C^{r}} \text{ has spectral gap property } \forall |t| \leq t_{1} \Big\}$$
be, then
$X$ is dense in $C(M , \R)$, in the uniform topology. Moreover, $C^{r}(M , \R)\setminus X$ is nowhere dense in the uniform topology; in other words, it is not dense in any nonempty open subset $U \subset C(M , \R)$.
\end{maincorollary}

In statistical mechanics, the parameter $t$ is often interpreted as the inverse of the temperature. In this way, we say that the potential $\phi$ has no low temperature phase transition if there is $t_{0} > 0$ such that $(-\infty , -t_{0}) \cup (t_{0} , +\infty) \ni t \mapsto P_{top}(f , t\phi)$ is analytic. Analogously, we say that the potential $\phi$ has no high temperature phase transition if there is $t_{1} > 0$ such that $(-t_{1} , t_{1}) \ni t \mapsto P_{top}(f , t\phi)$ is analytic. Thus, it follows from the previous Corollary: 

\begin{maincorollary}\label{maincorB}
    Let 
    $$X:= \Big\{\phi \in C^{r}(M , \R) : \phi \text{ has no low and high-temperature phase transition }\Big\}$$ be, then $X$ is dense in $C(M , \R)$, in the uniform topology. Moreover; $C^{r}(M , \R)\setminus X$
         is nowhere dense in the uniform topology.
\end{maincorollary}

For an intermittent skew-product in $TM2$ or $TM3$, defined previously, we obtain information on phase transitions with respect to any smooth potential:

\begin{maintheorem}\label{maintheoB}
    Suppose that $F \in TM2$ or $TM3$ and $\phi \in C^{r}(\T^{d} \times \Sc^{1}, \R)$. Then, there exists $A \subset \R$ open and dense subset such that $A \ni t \mapsto P_{top}(F , t\phi)$ is analytic.
\end{maintheorem}

For an intermittent skew-product in $TM1$, defined previously, we improve Theorem \ref{theorA}. In fact, we obtain the density of the spectral gap property for all parameters $t \in \R$, and the non-existence of phase transitions:

\begin{maintheorem}\label{maintheoC}
    Suppose that $F \in TM1$. 
    %and there exists $j \in \{1 , \ldots, k-1\}$, as in the Example \ref{examp}, such that $|f'_{x}(\alpha_{j}| > 1$ and $|f'_{x}(\alpha_{j+1})| > 1$ for all $x \in \mathbb{T}^{d}$. 
    Then, there exists $\tilde{\mathcal{H}} \subset C^{r}(\T^{d} \times \Sc^{1} , \R) $ dense subset in $C(\T^{d} \times \Sc^{1} , \R)$, such that: if $\phi \in \tilde{\mathcal{H}}$  then

(i) $f$ has no phase transition with respect to $\phi;$

    (ii)  $\mathcal{L}_{f,t\phi}$ has spectral gap property, acting on $C^{\alpha}(M ,\C)$, for all $ t \in \R$ and $0 < \alpha < 1$.
\end{maintheorem}

Note that, as well as the Corollaries \ref{maincorA} and \ref{maincorB},  the subset of the smooth potential that does not satisfy item (i) or (ii) of the previous theorem is nowhere dense in the uniform topology.

%\begin{maintheorem}\label{maintheoC}
 %   Suppose that $F \in TM1$ or $TM2$. Then, there exists $\tilde{\mathcal{H}} \subset C(M , \R)$ open and dense subset such that: if $\phi \in \tilde{\mathcal{H}} \cap C^{r}(M , \R)$ then $\mathcal{L}_{f,t\phi}$ has spectral gap property, acting on $C^{r}(M ,\C)$, for all $ t \in \R$ and $\R \ni t\mapsto P_{top}(F , t\phi)$ is strictly convex. In particular:

  %  (i) $\ds \Big\{\phi \in C^{r}(M , \R) : \mathcal{L}_{f,t\phi|C^{r}} \text{ has spectral gap property, for all } t \in \R\Big\} \cap \Big\{\phi \in C^{r}(M , \R) : f \text{ has no phase transition with respect to } \phi \Big\}$ is dense, in the uniform topology.

   % (ii) $\ds \Big\{\phi \in C^{r}(M , \R) : f \text{ has phase transition with respect to } \phi \Big\} \cup \Big\{\phi \in C^{r}(M , \R) : \mathcal{L}_{f,t\phi|C^{r}} \text{ has no spectral gap property for some } t \in \R  \Big\}$ is nowhere dense, in the uniform topology.
%\end{maintheorem}

\section{Preliminaries}\label{sec:preli}

In this section, we provide some definitions and preliminary results.

\subsection{Ergodic Theory}

Given a continuous transformation $T : X\rightarrow X$, acting on a compact metric space, we denote $\mathcal{M}_{1}(T) := \{\mu : \mu \text{ is a $T-$invariant probability} \}$ and  $\mathcal{M}_{e}(T) := \{\mu\in\mathcal{M}_{1}(T) : \mu \text{ is  $T-$ergodic} \}$.

In our context, we already know that $f : M \rightarrow M$ is topologically conjugated to a uniformly
expanding dynamic, in particular:

\begin{itemize}
    \item $f$ has the periodic specification and shadowing property;
    \item $\mathcal{M}_{1}(f) \ni \mu \mapsto h_{\mu}(f)$ is upper semicontinuous;
    \item $f$ has an equilibrium state with respect to every continuous potential,
    \item $\{f^{-n}(x) : n \geq 0\}$ is dense in $M$ for all $x \in M$, 
    \item $f$ admits a generating partition by domains of injectivity.
    \end{itemize}

For more details, see e.g. \cite{OV16}.

\subsection{Lyapunov exponents}

The Lyapunov exponents translate the asymptotic rates of expansion and contraction of a smooth dynamical system; these are defined via the Oseledets multiplicative ergodic theorem:

 %We say that $\lambda$ is a \textit{Lyapunov exponent} for $f : M \rightarrow M$, a $C^1$ map, if there exists a point $x$ and a vector  $v \in T_xM$ such that
%$$ \lambda=\lim_{n\to\infty} \dfrac{1}{n} \log \|Df^n_x(v)\|,$$

%We let $L(f)$ denote the set of all Lyapunov exponents for $f$.
%The Oseledets Ergodic Theorem states that 
For each ergodic measure $\mu \in \mathcal{M}_{e}(f)$
there exists constants $\lambda_{1}(f, \mu) > \dots > \lambda_{k}(f,\mu)$, and a filtration of $Df-$invariants subspaces  $T_xM = V_{1}(x) \supset \ldots \supset V_{k+1}(x) = \{0\}$,
such that 
$$\lim_{n \to \infty} \dfrac{1}{n} \log \|Df^n(x) v\|= \lambda_{i}(f , \mu) = \lambda_i$$
for $\mu$-almost every $x$ and every
vector $v \in V_{i}(x) \setminus V_{i+1}(x)$, $i = 1 \ldots, k$. For non-ergodic measures, the number $k$, the constants $\lambda_j$, and the tangent bundle decomposition, which depends on $x$, may depend on the ergodic component. The constants $\lambda_
j$ are called the Lyapunov exponents associated to the measure $\mu$ (for more details, see e.g. \cite{Rue79}). 
The maximum and minimum Lyapunov exponents $\lambda_{max}(\mu)=\lambda_1(\mu)$  and $\lambda_{min}(\mu)=\lambda_k(\mu)$ can be calculated by the norm and conorm of the derivative:
$$\lim_{n \to \infty} \dfrac{1}{n} \log \|Df^n(x)\|= \lambda_{max}(\mu)$$
$$\lim_{n \to \infty} \dfrac{1}{n} \log \|(Df^n(x))^{-1}\|^{-1}= \lambda_{min}(\mu)$$

%Now for an ergodic measure $\mu$, let $m_i:=m^i_x=\dim E^i_x-\dim E^{i+1}_x$ be the multiplicity for $i=1,\dots,k$   and $\mu-$almost every $x$. Define the sum of positive Lyapunov exponents, considering the multiplicity
%$$\lambda_+(\mu)=\sum_{\lambda_j>0} m_j \lambda_j.$$

%We will need the following well-known relation between entropy and positive Lyapunov exponents:

%\begin{theorem}[Margulis-Ruelle inequality, \cite{Rue78}]
%Let $f : M \rightarrow M$ be a $C^{1}$-local diffeomorphism that preserves an $f-$invariant and ergodic probability $\mu$. Then $$h_\mu(f)\leq  \lambda_+(\mu)$$
%\end{theorem}

\subsection{Transfer Operator}

All the necessary machinery  that relates the spectral properties of the transfer operator and the regularity of the pressure function used in \cite{BC21} and 
\cite{BF23} to prove the results for circle maps actually works for high-dimensional maps. Now we present such results, considering the same proofs as in \cite{BC21} or \cite{BF23}.

The first lemma states that if the transfer operator has a peripheral eigenvalue, then the spectral radius is a simple eigenvalue, and the associated eigenfunction is bounded away from zero.

\begin{lemma}\label{Lemaxi}\cite[Lemma 5.4]{BC21}
Let $f:M \rightarrow M$ be a $C^{r}-$local diffeomorphism on a compact and connected manifold $M$, such that $\{f^{-n}(x) : n \geq 0\}$ is dense in $M$ for all $x \in M$, and let $\phi \in C^{r}(M , \R)$ be a real function. If $\Lo_{f,\phi}\varphi=\lambda\varphi$
with $|\lambda|=\rho(\Lo_{f,\phi}|_{C^{r}})$ and $\varphi \in C^{r}(M , \R)\setminus\{0\}$, then
$\Lo_{f,\phi}|\varphi|=\rho(\Lo_{f,\phi}|_{C^{r}})|\varphi|.$ Furthermore, $\rho(\Lo_{f,\phi}|_{C^0})=\rho(\Lo_{f,\phi}|_{C^{r}})$, $\varphi$ is bounded away from zero and $\dim\ker(\Lo_{f,\phi}|_{C^{r}} - \lambda I)=1$.
\end{lemma}

With this bounded eigenfunction, we construct an invariant probability with full support.

\begin{corollary}\label{medida}\cite[Lemma 5.4 and Corollary 5.5]{BC21}
Let $f:M \rightarrow M$ be a $C^{r}-$local diffeomorphism on a compact and connected manifold $M$, such that $\{f^{-n}(x) : n \geq 0\}$ is dense in $M$ for all $x \in M$, and let $\phi \in C^{r}(M , \R)$ be a real function. If $\mathcal{L}_{f, \phi}|_{C^{r}}$ has the spectral gap property, then there exists a unique probability $\nu_{\phi}$ on $M$ such that $ (\mathcal{L}_{f, \phi}|_{C^{r}})^{\ast}\nu_{\phi} = \rho(\mathcal{L}_{f, \phi}|_{C^{r}})\nu_{\phi}$. Moreover, $supp(\nu_{\phi}) = M$.
\end{corollary}

\begin{remark}
Given $\Lo_{f,\phi}|_{C^{r}}$ with the spectral gap property,  by Lemma \ref{Lemaxi} there exists a unique $h_{\phi} \in C^{r}(M , \C)$ such that $\Lo_{f,\phi} h_{\phi}=\rho(\Lo_{f,\phi}) h_{\phi}$ and $\int h_{\phi} d\nu_{\phi} = 1$. Moreover, denote the $f-$invariant probability $h_{\phi}\nu_{\phi}$ by $\mu_{\phi}$. Note that $supp(\nu_{\phi}) = M$ and $h_{\phi} > 0$, thus, by Corollary \ref{medida}, $supp(\mu_{\phi}) = M.$
\end{remark}

The following result shows us that, in our context, the spectral gap property implies that the thermodynamic quantities are related to the spectral quantities.

\begin{lemma}\label{Lemapress}\label{GapAnalt}\cite[Lemma 4.13]{BF23}
 Let $f:M \rightarrow M$ be a $C^{r}-$local diffeomorphism on a compact and connected manifold $M$, such that $\{f^{-n}(x) : n \geq 0\}$ is dense in $M$ for all $x \in M$ and $f$ admits generating partition by domains of injectivity, and let $\phi \in C^{r}(M , \R)$ be a real function. If $\Lo_{f,\phi}|_{C^{r}}$ has the spectral gap property then:
\begin{enumerate}
    \item $P_{top}(f,\phi)=\log \rho(\Lo_{f,\phi}|_{C^{r}})$ and $\mu_{\phi}$ is the unique equilibrium state of $f$ with respect to $\phi$;
\item $\R \ni t\mapsto P(f,t\phi)$ is analytic in a neighbourhood of $1$.
\item  Suppose additionally that $\phi$ is not cohomologous to constant, in other words, there exists no $ c \in \R$ and $u \in C(M , \R)$ such that $\phi = c + u\circ f - u$. Then $\R \ni t \mapsto P_{top}(f , t\phi) $ is strictly convex in a neighborhood of $1$. 
\end{enumerate}

\end{lemma}

Another important fact is that, in our context, if the transfer operator is quasi-compact, then it has the spectral gap property. The proof is analogous to the respective result in \cite{BC21}. Recall that given $E$ a complex Banach space and $A:E \to E$ a bounded linear operator, define the essential spectral radius:
$$\rho_{ess}(A):=\inf\{r>0; \;\text{sp}(A)\setminus \overline{B(0,r)} \text{ contains only eigenvalues of finite multiplicity}\}.$$
\noindent Thus, we say that $A$ is quasi-compact if $\rho_{ess}(A) < \rho(A)$.

\begin{proposition}\label{LemaEss}\cite[Proposition 5.9]{BC21}
Let $f:M \rightarrow M$ be a $C^{r}-$local diffeomorphism  on a compact and connected manifold $M$, such that $\{f^{-n}(x) : n \geq 0\}$ is dense in $M$ for all $x \in M$, and let $\phi \in C^{r}(M , \R)$ be a real function. If $\Lo_{f,\phi}|_{C^r}$ is quasi-compact, then it has the spectral gap property.
\end{proposition}

We should also remember that it follows from \cite[Theorem 1.3]{W01}:

\begin{corollary}\label{pressaomenorqueraio}
$P_{top}(f,\phi) = \log \rho(\cL_{f,\phi}|_{C^{0}}).$
\end{corollary}

\begin{remark}
    
 Suppose that $\phi \in C^r(M , \R)$, so we know there exists a conformal measure $\nu_\phi$ associated with the spectral radius of $\Lo_{f,\phi}|_{C^{0}}$, thus $\nu_\phi \in \ker(\Lo^{\ast}_{f,\phi}|_{C^{r}} - \rho(\Lo^{\ast}_{f,\phi}|_{C^{0}})I)$ and $\rho(\Lo_{f,\phi}|_{C^{0}}) = \rho(\Lo^{\ast}_{f,\phi}|_{C^{0}})\in sp(\Lo^{\ast}_{f,\phi}|_{C^{r}}) = sp(\Lo_{f,\phi}|_{C^{r}})$. Therefore, $\rho(\Lo_{f,\phi}|_{C^{0}}) \leq \rho(\Lo_{f,\phi}|_{C^{r}})$.
\end{remark}

Later in this work, estimates of the essential spectral radius will be used to obtain quasi-compactness and consequently the spectral gap property. These estimates appear in \cite[Theorems 1 and 2]{CL97}:

\begin{theorem}\label{Lat}
Assume that $f:M\rightarrow M$ is any smooth covering and $\phi \in C^{r}(M, \R)$, then
$$ \rho_{ess}(\Lo_{f,\phi}|_{C^k}) \leq \exp\big[\sup_{\mu \in \mathcal{M}_{e}(f)}\{h_\mu(f)+\int\phi \text{d}\mu-k\lambda_{\min}(f , \mu)\}\big] \text{ and }$$
$$ \rho(\Lo_{f,\phi}|_{C^k}) \leq \exp\big[\sup_{\mu \in \mathcal{M}_{1}(f)}\{h_\mu(f)+\int\phi \text{d}\mu\}\big],$$
for $k=0,1,\dots, r.$ Where $\lambda_{\min}(f , \mu)$ denotes the smallest Lyapunov-Oseledec exponent of $\mu$.
\end{theorem}

\section{Proof of the main results}\label{sec:proof}

\subsection{Expanding on average potentials}
In order to construct our dense subsets of potentials, we consider \textit{expanding on average potentials}, which shall be defined and explored in this section.

\begin{remark}\label{nonneg}
     A map $f$ 
 in our setting has no negative Lyapunov exponents; this fact follows from Pesin's Theory for Endomorphisms. Indeed, suppose by contradiction that $\mu \in \mathcal{M}_{e}(f)$ has negative Lyapunov exponents. Hence, by the Stable Manifold Theorem (see \cite[Theorem V.6.4]{QXZ78}), there exists a stable manifold $\mathcal{W}^{s}(x) \supsetneqq \{x\}$. In particular, for all $\epsilon > 0$ there exists $y \neq x$ such that $d(f^{n}(x) , f^{n}(y)) < \epsilon, \forall n \geq 0$. This contradicts the fact that $f$ is expansive.
\end{remark}

Given this remark, we want to study measures with only positive Lyapunov exponents. So, following \cite{PV22}, we say that: 

\begin{definition}
A given $\mu \in \mathcal{M}_{1}(f)$ is an {\bf expanding measure} if all Lyapunov exponents of $\mu$ are positive, that is, 
$$\ds \lim_{n\to +\infty}\frac{1}{n}\log||Df^{n}(x) \cdot v|| > 0$$ for all $v \in T_{x}M\setminus \{0\}$ and $\mu-$a.e. $x \in M$. We will denote the expanding probability space by $\mathcal{E}(f)$. 
\end{definition}

In particular, if $\mu$ is expanding then 
$$\lim_{n \to +\infty} \log ||Df^{n}(x) \cdot v|| > 0$$ for all $||v || = 1$ and $\mu-$a.e. $x$.

Given the following \cite[Lemma A.4]{P11}:

\begin{lemma}
    Let $f : M \rightarrow M$ be a $C^{1}-$local diffeomorphism and let $\mu$ be an $f-$invariant and ergodic probability. If for $\mu$ almost every $x \in M$ we have 
    $$\lim_{n \to +\infty} \log ||Df^{n}(x) \cdot v|| > 0$$ for all $||v|| = 1$ then there exists $l \geq 1$ such that 
    $$\int \log ||Df^{l}(x)^{-1} ||^{-1}d\mu > 0. $$ 
    %$$\lim_{n \to +\infty} \frac{1}{n}\sum_{j = 1}^{n-1}\log ||Df^{l}(f^{lj}(x))^{-1} ||^{-1} > 0 $$ for $\mu$ almost every $x \in M$.
\end{lemma}

We can consider a weaker definition:

\begin{definition}
    We say that $\mu \in \mathcal{M}_{1}(f)$ is {\bf expanding on average} if there exists $l \geq 1$ such that 
    $$\int \log ||Df^{l}(x)^{-1} ||^{-1}d\mu > 0.$$
    We will denote the expanding on average probability space by $\mathcal{EA}(f)$.
\end{definition}

Thus, the previous Lemma implies that if $\mu$ is expanding, then $\mu$ is expanding on average. Furthermore, if $\mu \in \mathcal{M}_{e}(f)$ then both definitions are equivalent. Note that $\mathcal{EA}(f)$ is an open set, since it is a union of open sets.
We want to understand potentials whose equilibrium states are expanding (on average). 

\begin{definition}
 We say that $\phi \in C(M , \R)$ is an {\bf expanding potential } if 
$$P_{top}(f , \phi) > \sup\Big\{h_{\mu}(f) + \int\phi d\mu : \mu \in \mathcal{M}_{1}(f) \text{ and }\mu \notin \mathcal{E}(f)\Big\}.$$  

Similarly, we say that $\phi$ is an {\bf expanding on average potential} if 
$$P_{top}(f , \phi) > \sup\Big\{h_{\mu}(f) + \int\phi d\mu : \mu \in \mathcal{M}_{1}(f) \text{ and }\mu \notin \mathcal{EA}(f)\Big\}.$$ 
\end{definition}

Clearly, if a potential $\phi$ is expanding, then it is expanding on average, for any map $f$. We prove that in our context, both definitions are equivalent.

\begin{proposition}\label{wesp}
If $\phi \in C^r(M  ,\R)$ is an expanding on average potential, then $\Lo_{f,\phi}|_{C^{r}}$ has the spectral gap property. In particular, $\phi$ is an expanding potential.
\end{proposition}
\begin{proof}
Suppose $\Lo_{f,\phi}|_{C^{r}}$ does not have the spectral gap property; it follows from Proposition \ref{LemaEss} that it is not quasi-compact and the spectral and essential radius are equal. By Corollary \ref{pressaomenorqueraio} and Theorem \ref{Lat} we have
$$ P_{top}(f,\phi) = \log \rho(\cL_{f,\phi}|_{C^{0}}) \leq \log \rho(\cL_{f,\phi}|_{C^{r}}) = \log \rho_{ess}(\cL_{f,\phi}|_{C^{r}})\leq 
$$
$$
\sup_{\mu \in \mathcal{M}_{e}(f)}\left\{ h_\mu(f) + \int \phi d\mu - r\lambda_{\min}(\mu)\right\}.$$
Take a sequence $\mu_{n} \in \mathcal{M}_{e}(f)$ such that $ h_{\mu_{n}}(f) + \int \phi d\mu_{n} - r\lambda_{\min}(\mu_{n}) + \frac{1}{n} \geq \sup_{\mu \in \mathcal{M}_{e}(f)}\left\{ h_\mu(f) + \int \phi d\mu - r\lambda_{\min}(\mu)\right\},$ for all $n \in \N.$
In particular,
$$h_{\mu_{n}}(f) + \int \phi d\mu_{n} \leq P_{top}(f,\phi) \leq \sup_{\mu \in \mathcal{M}_{e}(f)}\left\{ h_\mu(f) + \int \phi d\mu - r\lambda_{\min}(\mu)\right\} \leq
$$
$$
h_{\mu_{n}}(f) + \int \phi d\mu_{n} - r\lambda_{\min}(\mu_{n}) + \frac{1}{n} \Rightarrow \lambda_{\min}(\mu_{n}) \leq \frac{1}{nr}.$$
It follows from Remark \ref{nonneg} that $\ds\lim_{n \to+\infty}\lambda_{\min}(\mu_{n}) = 0.$

Let $\ds \eta := \lim_{k \to +\infty}\mu_{n_{k}}$ be an accumulation point. For the sake of simplicity, reindex $\mu_{n_{k}}$ as $\mu_n$. Since the entropy function $\mu \mapsto h_{\mu}(f)$ is upper semicontinuous, we have:
$$
h_{\eta}(f) + \int \phi d\eta \geq \lim_{n\to +\infty}h_{\mu_{n}}(f) + \int \phi d\mu_{n}  = \lim_{n\to +\infty}h_{\mu_{n}}(f) + \int \phi d\mu_{n} - r\lambda_{\min}(\mu_{n}) + \frac{1}{n} \geq
$$
$$
P_{top}(f , \phi) 
\geq h_{\eta}(f) + \int \phi d\eta.
$$
Thus $\eta$ is an equilibrium state of $f$ with respect to $\phi$. Note that each $\mu_n$ has a generic point $x_n$, for which the space average above equals the time average. Hence, given $l \geq 1$ we have:
$$\int \log ||Df^{l}(x)^{-1} ||^{-1}d\eta = \lim_{n \to +\infty}\int \log ||Df^{l}(x)^{-1} ||^{-1}d\mu_{n} = $$
$$\lim_{n \to +\infty}
\lim_{m \to +\infty}\frac{1}{m}\sum_{j = 0}^{m-1}\log ||Df^{l}(f^{jl}(x_{n}))^{-1}||^{-1}\leq \lim_{n \to +\infty}\lim_{m \to +\infty}\frac{1}{m}\log||Df^{lm}(x_{n})^{-1}||^{-1} $$
$$=\lim_{n \to +\infty} l \lambda_{\min}(\mu_{n}) = 0 \Rightarrow \int \log ||Df^{l}(x)^{-1} ||^{-1}d\eta \leq 0.$$
In particular, $\eta$ is not expanding on average, and thus we conclude that 
$\phi$ is not expanding on average.

In particular, by Lemma \ref{Lemapress}, if $\phi \in C^r(M  ,\R)$ is an expanding on average potential, then $\phi$ has a unique equilibrium state, which is ergodic and expanding on average, thus expanding. Therefore, $\phi$ is an expanding potential.
\end{proof}

\subsection{Expanding on average potentials are common}

Given the Proposition \ref{wesp}, we want to study the denseness of the subset of expanding on average potentials in the uniform topology in order to prove Theorem \ref{theorA} 

Recall the Ergodic Optimization results and context of \cite{LT24}.
Let $T : X \rightarrow X$ be a continuous transformation on a compact metric space. Fix $c \in \Big[0 , h_{top}(T)\Big)$. Define $\Lambda_{c}:= \{\mu \in \mathcal{M}_{1}(T) : h_{\mu}(T) \geq c\}.$ 
Given $\phi \in C(X , \R)$, define:
$$
\mathcal{M}^{\Lambda_{c}}_{\max}(\phi) := \Big\{\mu \in \Lambda_{c} : \int \phi d\mu \geq \sup_{\nu \in \Lambda_{c}}\int \phi d\nu\Big\}.
$$
In our context, $T = f$ is transitive, it has the shadowing property and its associated entropy function $\mu \mapsto h_{\mu}(f)$ is upper semicontinuous. It follows from \cite[Theorem 3.1]{LT24}: 

\begin{theorem}
 If $L $ is residual in $\Lambda_{c}$ then $$D_{c}:= \Big\{\phi \in C(M , \R) : \mathcal{M}^{\Lambda_{c}}_{\max}(\phi) \subset L \text{ and } \#\mathcal{M}^{\Lambda_{c}}_{\max}(\phi) = 1\Big\}$$ is residual.
\end{theorem}

We call $\mu \in \mathcal{M}^{\Lambda_{c}}_{\max}(\phi)$ a \textbf{$c$-maximizing measure} and $\phi \in D_{c}$ a \textbf{uniquely $c$-maximized potential}. We now prove that the set of $c$-maximized expanding (on average) potentials is dense.

\begin{corollary}
Given $c  \in [0 , h_{top}(f))$, define 
$$\mathcal{E}_{c}:= \Big\{\phi \in C(M , \R) : \sup_{h_{\mu}(f) \geq c}\int \phi d\mu > \sup_{h_{\mu}(f) \geq c \text{ and } \mu \notin \mathcal{E}(f)}\int \phi  d\mu \Big\} \text{ and}$$
$$\tilde{\mathcal{E}}_{c}:= \Big\{\phi \in C(M , \R) : \sup_{h_{\mu}(f) \geq c}\int \phi d\mu > \sup_{h_{\mu}(f) \geq c \text{ and } \mu \notin \mathcal{EA}(f)}\int \phi  d\mu \Big\}.$$
Then $\mathcal{E}_{c}$ contains a residual subset, in particular $\mathcal{E}_{c}$ is a dense subset and $\tilde{\mathcal{E}}_{c}$ is an open and dense subset.
\end{corollary}
\begin{proof}
Take $L = \mathcal{EA}(f) \cap \Lambda_{c}$ from the previous theorem. 
    Since $\mathcal{EA}(f)$ open, then $L$ is an open subset in $\Lambda_{c}$. 
    
    Now, we must show that $L$ is a dense subset in $\Lambda_{c}$. Fix $\mu \in \mathcal{M}_{1}(f)$ with $h_{\mu}(f) \geq c$ and $\epsilon > 0$. Since $f$ is expansive, there exists $\mu_{0} \in \mathcal{M}_{e}(f)$ such that $h_{\mu_{0}}(f) = h_{top}(f).$ Taking $\mu_{s} := (1 - s)\mu + s\mu_{0}$, for a small enough $s > 0$ such that $\mu_s \in \mathcal{M}_{1}(f)$ with $d(\mu_s , \mu) < \frac{\epsilon}{4}$ and $h_{\mu_s}(f) > c.$
    Since $f$ has the specification property, it follows from \cite{PS05}, that $f$ is entropy-dense. In particular, there exists $\tilde{{\mu}} \in \mathcal{M}_{e}(f)$ such that $d(\tilde{{\mu}} , \mu_s)  < \frac{\epsilon}{4}$ and $h_{\tilde{{\mu}}}(f) > c.$ In particular, $d(\mu,\tilde{\mu})\leq d(\mu,\mu_s)+d(\mu_s,\tilde\mu)<\frac{\epsilon}{4}+\frac{\epsilon}{4} = \frac{\epsilon}{2}$
    
    Now, given $\tilde{\epsilon} > 0$ define:
    $$
    U_{\tilde{\epsilon}, n} := \Big\{x \in M :\frac{1}{n}\log ||Df^{n}(x)^{-1}||^{-1} \geq -\tilde{\epsilon} \Big\}.$$
   % and
   % $$
   % K_{n}:= \min_{x \in M}\frac{1}{n}\log ||Df^{n}(x)^{-1}||^{-1}.$$
    Hence
    $$
    \frac{1}{n}\int\log ||Df^{n}(x)^{-1}||^{-1} d\tilde{{\mu}}  =  
    $$
    $$
    \frac{1}{n}\int_{U_{\tilde{\epsilon},n}}\log ||Df^{n}(x)^{-1}||^{-1} d\tilde{{\mu}} \,+ \, \frac{1}{n}\int_{M\setminus U_{\tilde{\epsilon},n}}\log ||Df^{n}(x)^{-1}||^{-1} d\tilde{{\mu}}  \geq
    $$
    $$-\tilde{\epsilon}\cdot\tilde{{\mu}}(U_{\tilde{\epsilon},n}) \,+\, \inf_{x \in M}\log||Df(x)^{-1}||^{-1}\Big(1 - \tilde{{\mu}}(U_{\tilde{\epsilon},n})\Big).$$
    %$$-\tilde{\epsilon}\cdot \tilde{\tilde{\mu}}(U_{\tilde{\epsilon},n}) + K_{n}\Big(1 - \tilde{\tilde{\mu}}(U_{\tilde{\epsilon},n})\Big) \geq 
    %$$
    
    Recall that $\tilde{{\mu}}$ has no negative Lyapunov exponents, by Remark \ref{nonneg}. In particular, $\ds \lim_{n \to +\infty}\tilde{{\mu}}(U_{n , \tilde{\epsilon}}) = 1$ and the expression above tends to $-\tilde{\epsilon}$. Fixing $l$ the period of the repeller periodic point $p$, there exists $n_{\tilde{\epsilon}} > 0$ big enough such that
    $$\frac{1}{n_{\tilde{\epsilon}}l}\int\log ||Df^{n_{\tilde{\epsilon}}l}(x)^{-1}||^{-1} d\tilde{{\mu}} \geq -2\tilde{\epsilon}.$$
    
   Furthermore, define the expanding periodic measure $\nu := \frac{1}{l}\sum_{j = 0}^{l-1}\delta_{f^{j}(p)}.$ In particular, $\int \log ||Df^{l}(x)^{-1}||^{-1} d \nu > 0$.  
    
Let $\tilde{{\mu}}_{t} := (1-t)\tilde{{\mu}} + t\nu$. Then:

\begin{itemize}
    \item $\ds\frac{1}{n_{\tilde{\epsilon}}l}\int\log ||Df^{n_{\tilde{\epsilon}}l}(x)^{-1}||^{-1} d\tilde{{\mu}}_{t} \geq (1 - t)(-2\tilde{\epsilon}) + t \cdot \frac{\int \log ||Df^{l}(x)^{-1}||^{-1}d\nu}{l};$
    \item[] 
    \item $\ds d(\tilde{{\mu}} , \tilde{{\mu}}_t) = t d(\tilde{{\mu}},\nu) \leq t\Big( d(\tilde{\mu}, \mu) + d(\mu , \nu) \Big) \leq t \Big(\epsilon/2 + d(\mu , \nu)\Big);$
    \item[] 
    \item $\ds h_{\tilde{{\mu}}_{t}}(f) = (1 - t)h_{\tilde{{\mu}}}(f).$
\end{itemize}
Now, fix $t > 0$ close enough to $0$ such that $$t \Big(\epsilon/2 + d(\mu, \nu)\Big) < \frac{\epsilon}{2}\;\text{ and }\; (1 - t)h_{\tilde{{\mu}}}(f) \geq c.$$ Furthermore, fix $\tilde{\epsilon} > 0$ close enough to $0$ such that 
$$(1 - t)(-2\tilde{\epsilon}) + t \cdot \frac{\int \log ||Df^{l}(x)^{-1}||^{-1}d\nu}{l} > 0.$$
Finally, we have that $\tilde\mu_t \in L$ and is $\epsilon$-close to $\mu$:
\begin{itemize}
    \item $\ds\frac{1}{n_{\tilde{\epsilon}}l}\int\log ||Df^{n_{\tilde{\epsilon}}l}(x)^{-1}||^{-1} d\tilde{{\mu}}_{t}  > 0,$ in particular $\tilde{{\mu}}_{t}$ is expanding on average;
    \item[] 
    \item $\ds d(\tilde{{\mu}}_{t} , \mu) = d(\tilde{{\mu}}_{t} , \tilde{{\mu}}) + d(\tilde{{\mu}} , \mu) \leq t \Big( 2\epsilon + d(\mu , \nu)\Big) + d(\tilde{\mu} , \mu) < \frac{\epsilon}{2} + \frac{\epsilon}{2} = \epsilon;$
    \item[] 
    \item $\ds h_{\tilde{{\mu}}_{t}}(f) \geq c.$
\end{itemize}
We conclude that $L$ is a dense subset. \\

It follows from the previous theorem that $D_{c} := \Big\{\phi \in C(M , \R) : \mathcal{M}^{\Lambda_{c}}_{\max}(\phi) \subset L \text{ and } \#\mathcal{M}^{\Lambda_{c}}_{\max}(\phi) = 1\Big\}$ is residual. Recall that if $\mu \in \mathcal{M}_{e}(f)$ then $\mu $ is expanding on average if, and only if, $\mu$ is expanding. Therefore, $\mathcal{E}_{c}$ contains $D_c$, a residual subset. Moreover, we have that $\mathcal{E}_{c} \subset \tilde{\mathcal{E}}_{c}$ and since $\mathcal{EA}(f)$ is an open subset, then the set $\mathcal{M}_{1}(f) \setminus \mathcal{EA}(f)$ is closed. Finally, since $\mathcal{M}_{1}(f) \ni \mu \mapsto \int\phi d\mu $ is continuous and $\mathcal{M}_{1}(f) \ni \mu \mapsto h_{\mu}(f) $ is upper semicontinuous, we have that $\tilde{\mathcal{E}}_{c}$ is an open subset.

\end{proof}

Now we prove the desired density property.

\begin{proposition}\label{propexpd}
     The subset $\Big\{ \phi \in C(M , \R) : \phi \text{ is expanding }\Big\}$ is dense.
\end{proposition}
\begin{proof}
We will initially show that, for a dense subset of potentials $\phi$, the function $t \mapsto P_{top}(f , \phi)$ is analytic, strictly convex and $t\phi$ has a unique equilibrium state, for all $t \in \R$. Since $f$ is expansive, let $g : N \rightarrow N$ be an expanding map and $h : N \rightarrow M$ be a homeomorphism such that $f \circ h = h \circ g$. We know that 
    $$\tilde{C} := 
    \Big\{\phi : N \rightarrow \R : \phi \text{ is H\"older and not cohomologous to constant} \Big\}$$
    is dense. Furthermore, if $\phi \in \tilde{C}$ then $\R \ni t \mapsto P_{top}(g , t\phi)$ is analytic, strictly convex and $t\phi$ has a unique equilibrium state, for all $t \in \R$ (see e.g. \cite{PU10}). In particular, taking $C:= \{ \tilde{\phi} \circ h^{-1}\}_{\tilde{\phi} \in  \tilde{C}} $ we have that $C$ is dense and for all $\phi \in C$ the topological pressure function $\R \ni t \mapsto P_{top}(f , t\phi)$ is analytic, strictly convex and $t\phi$ has a unique equilibrium state, for all $t \in \R$.

    Given $\phi \in C(M , \R)$ and $\epsilon > 0$, take $\tilde{\phi} \in C$ such that $||\tilde{\phi} - \phi||_{\infty} < \frac{\epsilon}{17}$. For each $t\in \R$, let $\mu_{t}$ be the unique equilibrium state to $t\tilde{\phi}$. Furthermore, $\frac{dP_{top}(f , t\tilde{\phi})}{dt} = \int \phi d\mu_{t}$ by \cite{W92} and 
    $$
        h_{\mu_{t}}(f) = P_{top}(f,t\tilde{\phi}) - t\int\tilde{\phi} d\mu_t = P_{top}(f,t\tilde{\phi}) - t\frac{dP_{top}(f , t\tilde{\phi})}{dt} \Rightarrow 
        $$
        $$\frac{dh_{\mu_{t}}(f)}{dt} = -t\frac{d^{2}P_{top}(f , t\tilde{\phi})}{dt^{2}}. 
    $$
    
    Since $\R \ni t \mapsto P_{top}(f , t\tilde{\phi})$ is strictly convex, there exists $a > 0$ such that $\frac{d^{2}P_{top}(f , t\tilde{\phi})}{dt^{2}} \geq a$, for all $ t \in[\frac{1}{2} , 2]$.  Hence $\frac{dh_{\mu_{t}}(f)}{dt} \leq -\frac{a}{2} < 0$, for all  $ t \in[\frac{1}{2} , 2]$. We conclude that there exists $K > 0$ such that 
    $$
    |h_{\mu_{t_{1}}}(f) - h_{\mu_{t_{2}}}(f)| \geq K|t_{1} - t_{2}|, \, \forall t_1, t_2 \in[\frac{1}{2} , 2].
    $$
Since $\tilde{\mathcal{E}}_{h_{\mu_{1}}(f)}$ is an open and dense subset, by the previous corollary, take a sequence $\tilde{\tilde{\phi}}_{k} \in \tilde{\mathcal{E}}_{h_{\mu_{1}}(f)} \cap C$ such that $\lim ||\tilde{\tilde{\phi}}_{k} - \tilde{\phi}||_{\infty} = 0.$ For each $t \in \R$, let $\mu_{t,k}$ be the unique equilibrium state to $\tilde{\tilde{\phi}}_{k}$. Since $f$ is expansive, every accumulation point of $\mu_{t,k}$ is an equilibrium state to $\tilde{\phi}.$ Recall that $t\tilde{\phi}$ has a unique equilibrium state $\mu_t$, thus $\lim \mu_{t,k} = \mu_{t}$ and 
$$
\lim_{k \to +\infty}h_{\mu_{t,k}}(f) = \lim_{k \to +\infty}P_{top}(f,t\tilde{\tilde{\phi}}_{k}) - t\int\tilde{\tilde{\phi}}_{k} d\mu_{t,k} = P_{top}(f,t\tilde{\phi}) - t\int\tilde{\phi} d\mu_t
$$
Thereby, for $k$ large enough:
$$
    |h_{\mu_{t_{1},k}}(f) - h_{\mu_{t_{2},k}}(f)| \geq \frac{K}{2}|t_{1} - t_{2}|, \, \forall t_{1}, t_{2} \in[\frac{1}{2} , 2].
    $$
In particular, there exist $k_{0}$ large enough and $t_{k_{0}} \in \R$ such that $h_{\mu_{t_{k_{0}},k_{0}}}(f) = h_{\mu_{1}}(f)$ and $|t_{k_{0}} - 1|< \frac{\epsilon}{17}$.

Now we will show that $t_{k_{0}}\tilde{\tilde{\phi}}_{k_{0}}$ is an expanding potential. In fact, take $\nu \in \mathcal{M}_{1}(f)$ with $h_{\nu}(f) \geq h_{\mu_{1}}(f).$ Hence:
$$
h_{\mu_{1}}(f) + t_{k_{0}}\int \tilde{\tilde{\phi}}_{k_{0}} d\nu \leq h_{\nu}(f) + t_{k_{0}}\int \tilde{\tilde{\phi}}_{k_{0}}d\nu \leq h_{\mu_{t_{k_{0}}, k_{0}}}(f) + t_{k_{0}}\int \tilde{\tilde{\phi}}_{k_{0}}d\mu_{t_{k_{0}}, k_{0}} =
$$
$$
h_{\mu_{1}}(f) + t_{k_{0}}\int \tilde{\tilde{\phi}}_{k_{0}}d\mu_{t_{k_{0}}, k_{0}} \Rightarrow 
\int \tilde{\tilde{\phi}}_{k_{0}} d\nu \leq \int \tilde{\tilde{\phi}}_{k_{0}}d\mu_{t_{k_{0}}, k_{0}}.
$$
    Thus $\mu_{t_{k_{0}}, k_{0}}$ is $h_{\mu_{1}}(f)$-maximizing for $\tilde{\tilde{\phi}}_{k_{0}}$.  Since $\tilde{\tilde{\phi}}_{k_{0}} \in \tilde{\mathcal{E}}_{h_{\mu_{1}}(f)} \cap C,$ we have that $\mu_{t_{k_{0}}, k_{0}}$ is expanding and $t_{k_{0}}\tilde{\tilde{\phi}}_{k_{0}} $ is an expanding potential. 

This concludes the proof of the proposition.

\end{proof}

\subsection{Avoiding low-temperature phase transitions}

Let $$\tilde{\mathcal{E}}_{0}:=\Big\{\phi \in C(M , \R) :  \sup_{\mu \notin \mathcal{EA}(f)}\int \phi d\mu < \sup_{\mu \in \mathcal{EA}(f)}\int \phi d\mu\Big\}.$$
We already know that $\tilde{\mathcal{E}}_{0}$ is an open and dense subset.

The next lemma states that every potential in this set becomes expanding on average for large enough $t$.

\begin{lemma}\label{lemma5}
If $\phi \in \tilde{\mathcal{E}}_{0}$, then there exists a parameter $t_0>0$ such that $t\phi$ is expanding on average for all $t\geq t_0$. 
\end{lemma}
\begin{proof}
Suppose that the result does not hold, and take $t_n \to \infty$ with $t_n\phi$ not expanding on average. Let $\mu_n$ be an equilibrium state with respect to $t_n\phi$ such that $\int \log ||Df^{l}(x)^{-1} ||^{-1}d\mu_{n} \leq 0$, for all $l \geq 1$. Take an accumulation point $\eta = \lim\mu_{n_{k}}$.
In particular, $\eta$ is not expanding on average.
Furthermore, given $\mu \in \mathcal{M}_{1}(f)$:
$$h_\mu(f)+t_n\int \phi d\mu \leq h_{\mu_{n}}(f)+t_n\int\phi d\mu_{n}$$
$$\Rightarrow \dfrac{h_\mu(f)}{t_n} + \int \phi d\mu \leq \dfrac{h_{\mu_{n}}(f)}{t_n} + \int \phi d\mu_n $$
Taking the limit as $n\to \infty$
$$ \int \phi d\mu \leq \int \phi d\eta.$$
Therefore, $\phi \notin \tilde{\mathcal{E}}_{0}.$
\end{proof}

\begin{proposition}\label{proplt}
There exists $\tilde{D} \subset C(M , \R)$ an open and dense subset such that if $\phi \in \tilde{D} \cap C^{r}(M , \R)$ then there exists $t_{0} > 0$ such that $\mathcal{L}_{t\phi|C^{r}}$ has the spectral gap property $\forall |t| \geq t_{0}$ and $(-\infty , -t_{0}) \cup (t_{0} , +\infty) \ni t \mapsto P_{top}(f , t\phi)$ is analytic and strictly convex.
\end{proposition}
\begin{proof}
    It's enough to take $\tilde{D} := \tilde{\mathcal{E}}_{0}\cap (-\tilde{\mathcal{E}}_{0}),$ and apply Lemma \ref{lemma5}, Proposition \ref{wesp} and Lemma \ref{Lemapress}. Note that if $\phi \in \tilde{\mathcal{E}}_{0} $ then $\phi$ is not cohomologous to constant.
\end{proof}

\subsection{Avoiding high-temperature phase transitions}

In this section, we will prove that given $\phi \in C^{r}(M , \R)$ there exists $t_{1} > 0$ such that $\mathcal{L}_{t\phi|C^{r}}$ has the spectral gap property, for all $|t| \leq t_{1}.$ It follows from Lemma \ref{Lemapress} that $(-t_{1} , t_{1}) \ni t \mapsto P_{top}(f, t\phi)$ is analytic.

\begin{remark}
    Since $f$ is topologically conjugated to an expanding map then $f$ satisfies the periodic specification property, in other words, given $\epsilon >0$ there exists $N(\epsilon) > 0$ such that: if $x_{1}, \ldots, x_{k} \in M$ and $n_{1}, \ldots , n_{k} \in \N$ then there exists $y \in M$ with:
    \begin{itemize}
        \item $d(f^{i}(y) , f^{i}(x_{1}) < \epsilon,$ for all $i = 0 ,\ldots, n_{1};$
        \item $d(f^{n_{1}+ N(\epsilon)+ i}(y) , f^{i}(x_{2})) < \epsilon,$ for all $i = 0 ,\ldots, n_{2};$
        \item $d(f^{\sum_{j=1}^{l}n_{j}+ lN(\epsilon)+ i}(y) , f^{i}(x_{l+1})) < \epsilon,$ for all $l < k$ and $i = 0 ,\ldots, n_{l+1};$
        \item $f^{\sum_{j=1}^{k}n_{j}+ kN(\epsilon)}(y) = y.$
    \end{itemize}

\end{remark}

\begin{proposition}\label{propht}
There exists $\epsilon > 0$ such that: if $\phi \in C(M , \R)$ satisfies $||\phi||_{\infty} < \epsilon$  then $\phi$ is expanding on average. In particular, if $\phi \in C^{r}(M , \R)$ satisfies $||\phi||_{\infty} < \epsilon$ then $\mathcal{L}_{\phi|C^{r}}$ has the spectral gap property.
\end{proposition}
\begin{proof}
Since be expanding on average potential is an open property and by Proposition \ref{wesp}, it is enough to prove that $\phi \equiv 0$ is an expanding potential. Let $p \in M$ be a repeller periodic point with $f^{k}(p) = p.$ We are going to show that the unique maximum entropy measure of $g := f^{k}$ is expanding. In particular, the unique maximum entropy measure of $f$ will be expanding. Note that $g$ is an expansive $C^{r}-$local diffeomorphism with the repelling fixed point $p$. Therefore, we can assume that $k = 1.$

Let $x_{1}, \ldots, x_{deg(g)} \in M$ be such that $x_{1} = p$ and $g(x_{j}) = p$, for all $j = 1, \ldots deg(g).$ Take $\epsilon > 0$ small enough such that: 
\begin{itemize}
    \item $B(x_{i}, \epsilon) \cap B(x_{j}, \epsilon) = \emptyset,$ for all $i \neq j;$
    \item there exists $ a > 0$ such that $\log ||Dg(x)^{-1}||^{-1} \geq a > 0,$ for all $x \in B(x_{1}, \epsilon)$.
\end{itemize}

Fix $N(\epsilon) >0$, given by periodic specification property for $g$. Take $n_{0} \in \N$ such that:
$$
\frac{-N(\epsilon) \inf _{x \in M}\log ||Dg(x)^{-1}||^{-1}}{n_{0}a } < 1.
$$
Fix $\max\{0 , \frac{-N(\epsilon) \inf _{x \in M}\log ||Dg(x)^{-1}||^{-1}}{n_{0}a - \inf _{x \in M}\log ||Dg(x)^{-1}||^{-1}}\} < \theta < 1.$
Thus; for all $k \geq 1$ and choices $x_{j_{1}}, \ldots, x_{j_{k}},$ with $j_{i} \in \{1 , \ldots, deg(g)\}$ and $\frac{\#\{j_{i} : j_{i} = 1\}}{k} \geq \theta$, there exists $y_{j_{1}, \ldots, j_{k}} \in M$ such that 
 \begin{itemize}
        \item $d(g^{i}(y_{j_{1}, \ldots, j_{k}}) , g^{i}(x_{j_{1}}) < \epsilon,$ for all $i = 0 ,\ldots, l_{1};$
        \item $d(g^{l_{1}+ N(\epsilon)+ i}(y_{j_{1}, \ldots, j_{k}}) , g^{i}(x_{j_{2}})) < \epsilon,$ for all $i = 0 ,\ldots, l_{2};$
        \item $d(g^{\sum_{r=1}^{m}l_{r}+ mN(\epsilon)+ i}(y) , f^{i}(x_{m+1})) < \epsilon,$ for all $m < k$ and $i = 0 ,\ldots, l_{m+1};$
        \item $g^{\sum_{r=1}^{k}l_{r}+ kN(\epsilon)}(y_{j_{1}, \ldots, j_{k}}) = y_{j_{1}, \ldots, j_{k}};$
    \end{itemize}
where $l_{r} = n_{0},$ if $j_{r} = 1,$ otherwise $l_{r} = 0.$

By construction;
$$
\frac{\sum_{r = 0}^{kN(\epsilon) + \sum_{i = 1}^{k}l_{i}-1}\log ||Dg(g^{r}(y_{j_{1}, \ldots, j_{k}})^{-1}||^{-1}}{kN(\epsilon) + \sum_{i = 1}^{k}l_{i}} \geq \frac{\theta kn_{0} a + kN(\epsilon)\inf_{x \in M}\log ||Dg(x)^{-1}||^{-1} }{kN(\epsilon) + \sum_{i = 1}^{k}l_{i}} 
$$
$$
\geq \frac{\theta kn_{0} a + kN(\epsilon)\inf_{x \in M}\log ||Dg(x)^{-1}||^{-1} }{kN(\epsilon) + kn_{0}} = 
$$
$$
\frac{\theta n_{0} a + N(\epsilon)\inf_{x \in M}\log ||Dg(x)^{-1}||^{-1} }{N(\epsilon) + n_{0}}  > 0;
$$
in particular $y_{j_{1}, \ldots, j_{k}}$ is a repeller periodic orbit. Note that $(j_{1}, \ldots, j_{k}) \mapsto y_{j_{1}, \ldots, j_{k}}$ is injective.

Since $g$ is topologically conjugate to an expanding map, we already know that $g$ has a unique maximum entropy measure and $h_{top}(g) = \log deg(g).$ 
Furthermore, 
$$\mu = \lim_{n \to +\infty}\frac{\sum_{g^{n}(x) = x}\delta_{x}}{\#\{x \in M : g^{n}(x) = x\}} = \lim_{n \to +\infty}\frac{1}{n}\frac{\sum_{i=0}^{n-1}\sum_{g^{n}(x) = x}\delta_{g^{i}(x)}}{\#\{x \in M : g^{n}(x) = x\}}$$ is the unique maximum entropy measure of $g$ (see \cite{Bo74}).
Define 
$$\lambda := \frac{\theta n_{0} a + N(\epsilon)\inf_{x \in M}\log ||Dg(x)^{-1}||^{-1} }{N(\epsilon) + n_{0}}$$ and 
$$E_{n} := \{x \in M : g^{n}(x) = x \text{ and } \frac{1}{n}\sum_{i=0}^{n-1}
\log||Dg(g^{i}(x))^{-1}||^{-1} \geq \lambda\}.$$ By previous construction;
$\limsup\frac{\#E_{n}}{\#\{x \in M : g^{n}(x) = x\}} >0.$ Therefore,
$$
\int\log||Dg(x)^{-1}||^{-1}d\mu = \limsup \frac{1}{n}\frac{\sum_{i=0}^{n-1}\sum_{g^{n}(x) = x}\delta_{g^{i}(x)}}{\#\{x \in M : g^{n}(x) = x\}} \geq 
$$
$$
\lambda \cdot \limsup\frac{\#E_{n}}{\#\{x \in M : g^{n}(x) = x\}} > 0.
$$
In particular, $\mu$ is expanding.
\end{proof}

\subsection{Proof of Theorem \ref{theorA}}

Define 
$$\mathcal{H} := \{\phi \in C(M , \R): \phi \text{ is expanding on average } \} \cap \tilde{D},$$ where $\tilde{D}$ is given by Proposition \ref{proplt}. Note that $$\{\phi \in C(M , \R): \phi \text{ is expanding on average } \}$$ is open and dense by Proposition \ref{propexpd}, and $\tilde{D}$ is open and dense by Proposition \ref{proplt}. Thus, applying Proposition \ref{wesp}, Proposition \ref{proplt} and Proposition \ref{propht}, we conclude the proof of Theorem \ref{theorA}.

\subsection{Proof of Theorem \ref{maintheoB}} 

In this section $F \in TM1, TM2 \text{ or } TM3$, according to Section \ref{subsec:set}.

Given a probability $\mu \in \mathcal{M}_{e}(F)$ the Lyapunov exponents will be $$\lambda_{1}(\mu), \ldots, \lambda_{d}(\mu), \lambda^{c}(\mu),$$ where $\lambda_{i}(\mu)$ are the Lyapunov exponents of $\pi^{\ast}\mu \in \mathcal{M}_{e}(g)$ for $\pi(x , y) := x$ and $\lambda^{c}(\mu) := \int \log |f_{x}'(y)| d\mu(x,y)$. Note that if $g$ is intermittent, in the definition of $F$, then $d = 1$ and $\lambda_{1}(\mu) = \int \log |g'|d(\pi^{\ast}\mu).$

In this moment, we will suppose that $F \in TM2 \text{ or } TM3$. Suppose that $\mu \in \mathcal{M}_{e}(F)$ is not expanding. We will show that $\mu = \nu \times \delta_{q_{j}}$ or $\delta_{p_{j}}\times \eta,$ where $\nu \in \mathcal{M}_{e}(g)$, $\eta\in\mathcal{M}_{e}(f_{p_{j}})$ and $p_{1}, \ldots, p_{l}$ are the fixed points of $g$ and $q_1, \dots, q_{k}$ fixed points of $f_x$.
In fact, if $\lambda^{c}(\mu) = \int \log |f_{x}'(y)| d\mu(x,y) = 0$, since $\log |f_{x}'(y)| \geq 0$ 
for all $(x , y)\in \mathbb{T}^{d} \times \mathbb{S}^{1},$ then $\log |f_{x}'(y)| = 0$ for $\mu-$a.e. $(x , y).$ Note that if $|f_{x}'(y)| = 1$ then 
$y = q_{j},$ for some $j = 1 , \ldots, k.$ Thus $\mu = \pi^{\ast}\mu \times \delta_{q_{j}}.$ On the other hand, if $\lambda_{1}(\mu) = \int \log |g'|d(\pi^{\ast}\mu) = 0$, analogously to previous case, then $\mu(\{x_{j}\} \times \mathbb{S}^{1}) = 1.$ Thus $\mu = \delta_{p_{j}}\times \eta,$ for $\eta\in\mathcal{M}_{e}(f_{p_{j}}).$

\

\begin{proof}[Proof of Theorem \ref{maintheoB}]
    Let $\phi \in C^{r}(\mathbb{T}^{d} \times \mathbb{S}^{1} , \R)$. Define the following subsets:
$$
A_{1} := \Big\{t \in \R : \mathcal{L}_{t\phi|C^{r}} \text{ has the spectral gap property }\Big\}, 
$$
$$
A_{2,j}:=\Big\{ t \in \R : t\phi \text{ admits an equilibrium state } \mu = \nu \times \delta_{q_{j}}\Big\},$$
$$
A_{3,j}:=\Big\{ t \in \R : t\phi \text{ admits an equilibrium state } \mu = \delta_{p_{j}}\times \eta \Big\}.$$

\

\emph{Claim: $\R = A_{1} \cup (\bigcup_{j=1, \ldots k}A_{2,j}) \cup (\bigcup_{j = 1, \ldots, l}A_{3,j}).$}

\

\noindent In fact; let $t \in \R$ be such that $\mathcal{L}_{t\phi|C^{r}}$  has no the spectral gap property. Following the proof of the Proposition \ref{wesp}, we find $\eta$ an equilibrium state with respect to $t\phi$ such that $\eta = \lim \mu_{n}$ with $\mu_{n} \in \mathcal{M}_{e}(F)$ and $\lim \lambda_{\min}(\mu_{n}) = 0$. In our context, $\lambda_{\min}(\mu_{n}) = \int \log |f_{x}'(y)| d\mu_{n}(x,y)$ or $\int \log |g'|d(\pi^{\ast}\mu_{n})$. Hence $\int \log |f_{x}'(y)| d\eta(x,y) = 0$ or $\int \log |g'|d(\pi^{\ast}\eta) = 0$. Since $\log |f_{x}'(y)| \geq 0$ and $\log |g'(x)| \geq 0$ for all $(x , y) \in \mathbb{T}^{d} \times \mathbb{S}^{1}$, applying the Ergodic Decomposition Theorem (see e.g. \cite{V14}), there exists $\mu \in \mathcal{M}_{e}(F)$ such that $\mu$ is an equilibrium state with respect to $t\phi$ and $\lambda_{\min}(\mu) = 0$. Thus, follows from previous discussion that $t \in (\bigcup_{j=1, \ldots k-1}A_{2,j}) \cup (\bigcup_{j = 1, \ldots, l-1}A_{3,j}).$\\

Moreover, note that:

\begin{itemize}
    \item $A_{1}$ is open and $A_{1} \ni t \mapsto P_{top}(f , t\phi)$ is analytic, by Lemma \ref{Lemapress}.
    \item[] 
    \item If $t \in A_{2,j}$ take $\nu \times \delta_{q_{j}} $ an equilibrium state with respect to $t\phi$. Thus, for all $\xi \in \mathcal{M}_{1}(g)$ we have:
    $$
    h_{\xi \times \delta_{q_{j}}}(F) + t \int \phi d (\xi \times \delta_{q_{j}}) \leq h_{\nu \times \delta_{q_{j}}}(F) + t \int \phi d(\nu \times \delta_{q_{j}}) \Rightarrow
    $$
    $$
    h_{\xi }(g) + t \int \phi(x , q_{j}) d \xi(x) \leq h_{\nu }(g) + t \int \phi(x , q_{j}) d\nu(x).
    $$
    We conclude that $\nu$ is an equilibrium state of $g$ with respect to $t\phi(\cdot , q_{j})$ and $P_{top}(F , t\phi) = P_{top}\big(g , t\phi( \cdot, q_{j})\big).$
    \item[]
    \item On the other hand, if $t \in A_{3,j}$ take $\delta_{p_{j}}\times \eta $ an equilibrium state with respect to $t\phi$. Analogously the previous item, for all $\xi \in \mathcal{M}_{1}(f_{p_{j}})$ we have:
    $$
    h_{\xi}(f_{p_{j}}) + t \int \phi(p_{j} , y) d \xi(y) = h_{\delta_{p_{j}}\times \xi}(F) + t \int \phi d (\delta_{p_{j}}\times \xi) \leq
    $$
    $$
    h_{\delta_{p_{j}}\times \eta}(F) + t \int \phi d(\delta_{p_{j}}\times \eta) = h_{ \eta}(f_{p_{j}}) + t \int \phi(p_{j} , y) d \eta(y). 
    $$
    We conclude that $\eta$ is an equilibrium state of $f_{p_{j}}$ with respect to $t\phi(p_{j}, \cdot)$ and $P_{top}(F , t\phi) = P_{top}\big(f_{p_{j}} , t\phi( x_{j} , \cdot )\big).$
\end{itemize}

Therefore, take $A:= A_{1} \cup \big(\bigcup_{j=1, \ldots k-1}int (A_{2,j})\big) \cup (\bigcup_{j = 1, \ldots, l-1}int (A_{3,j})\big), $ where $int(B)$ means the interior of subset $B$. Note that:

\begin{itemize}
    \item $A$ is an open subset.
    \item[]
    \item The function $int (A_{2,j}) \in t \mapsto P_{top}(F , t\phi) = P_{top}(g , t\phi(\cdot , q_{j})$ is analytic if $g$ is expanding, otherwise follows from \cite{BF23} that the function is analytic except at most two parameters $t$.
    \item[]
    \item The function $int (A_{3,j}) \in t \mapsto P_{top}(F , t\phi) = P_{top}(f_{p_{j}} , t\phi(p_{j} , \cdot)$ is analytic except at most two parameters $t$, by \cite{BF23}.
\end{itemize}

\

\emph{Claim: $A$ is dense in $\R$.}

\

\noindent In fact; let $\tilde{t} \notin A$ be and $\epsilon > 0$. Then there exist $|\tilde{\tilde{t}} - \tilde{t}| < \epsilon $ and $j$ such that 
$$P_{top}(g  , \tilde{\tilde{t}}\phi(\cdot, q_{j})) > \max\Big\{ \max_{i \neq j}P_{top}(g , \tilde{\tilde{t}}\phi(\cdot , q_{i})) , \max_{i}P_{top}(f_{p_{i}} , \tilde{\tilde{t}}\phi(p_{i}, \cdot))\Big\}$$
or
$$P_{top}(f_{p_{j}}  , \tilde{\tilde{t}}\phi(p_{j},\cdot)) > \max\Big\{ \max_{i \neq j}P_{top}(f_{p_{i}}  , \tilde{\tilde{t}}\phi(p_{i},\cdot)), \max_{i}P_{top}(g , \tilde{\tilde{t}}\phi(\cdot, q_{i}))\Big\}.$$
Since the previous inequalities are preserved for $t$ close enough to $\tilde{\tilde{t}}$, if $\tilde{\tilde{t}} \notin A$ then there exits a sequence $t_{n} \mapsto \tilde{\tilde{t}}$ with $t_{n} \in A.$ Thus $A$ is dense.

We conclude that $A$ is dense and the function $A \in t \mapsto P_{top}(F , t\phi)$ is analytic except at most in a finite number of parameters $t$.
\end{proof}

\

\subsection{Proof of Theorem C}

To prove Theorem \ref{maintheoC}, we need to remember the context and results in \cite{Kl20}, guaranteeing a class of potentials whose associated transfer operator has the spectral gap property. We can assume without loss of generality that $\T^{d}\times \Sc^{1}$ is endow with the distance $d((x_{1} , y_{1}) , (x_{2}, y_{2})) := d_{\T^d}(x_{1}, x_{2}) + \alpha \cdot |y_{1}- y_{2}|,$ where $\alpha > 0$ will be chosen to be small enough.

Fix $T : \Omega \rightarrow \Omega$ a local diffeomorphism on a compact and connected manifold, with topological degree $k$. Take the uniform backward random walk $\textbf{M} = (m_{x})_{x \in \Omega}$  defined by $\ds m_{x} = \frac{1}{k}\sum_{y \in T^{-1}(x)}\delta_{y}.$

\begin{definition}\label{def:transker}
    A continuous map $c : [0 , +\infty) \rightarrow [0 , +\infty)$ will be called \emph{contraction function} if $c(0) = 0$ and $c(r) < r$ for all $r >0.$ We say that the uniform backward random walk $\textbf{M}$ is \emph{weakly contracting} if there exists a contraction function $c$ and a real number $\lambda > 1$ such that, for all $x,y \in \Omega$, there exist permutations $\eta, \sigma$ of $\{1, \ldots, k\}$ such that:
    \begin{itemize}
        \item[(i)] for all $j \in \{1, \ldots, k\},$ we have $d(x^{\eta(j)}, y^{\sigma(j)}) \leq c\big( d(x , y)\big),$ where $x^{\eta(j)} \in T^{-1}(x)$ and $y^{\sigma(j)} \in T^{-1}(y)$;
        \item[(ii)] $d(x^{\eta(k)}, y^{\sigma(k)}) \leq \frac{d(x , y))}{\lambda}.$
    \end{itemize}
    \end{definition}

In our context, we want to prove that almost expanding skew-products in $TM1$ produce weakly contracting random walks.

\begin{proposition}
    For $F:\T^d\times\Sc^1 \to \T^d\times\Sc^1 \in TM1$, the associated uniform backward random walk {\bf M} is weakly contracting. 
\end{proposition}
\begin{proof}
On the one hand; note that fixed $\epsilon > 0$ there exists  $\lambda_{\epsilon} > 1$ such that $\inf_{||v|| = 1}||DF(x , y) \cdot v|| \geq \lambda_{\epsilon}$, if 
$\ds\inf_{j = 1, \ldots, k }|y - q_{j,x}| \geq \epsilon$; moreover, given $(x, y) \in \T^{d}\times \Sc^{1}$ there exists $(x',y') \in F^{-1}(x , y)$ such that $\ds\inf_{j = 1, \ldots, k }|y' - q_{j,x'}| \geq \epsilon$.

On other one hand; let's remember that since $F$ is $C^{1}$ there exists $\sigma >0$ such that $|f_{x_{1}}(y) - f_{x_{2}}(y)| \leq \sigma d_{\T^{d}}(x_{1},x_{2}),$ for all $x_{1} , x_{2} \in \T^{d}$ and $y \in \Sc^{1}$. Moreover, since $g$ is expanding there exists $\lambda > 1$ such that: for all $x_{1} , x_{2} \in \T^{d}$ we have $g^{-1}(x_{1}) = \{w_{1,1}, \ldots, w_{1,deg(g)}\}$ and $w_{2} \in g^{-1}(x_{2}) = \{w_{2,1}, \ldots, w_{2,deg(g)}\}$ with $d(w_{1,j} , w_{2,j}) \leq \frac{d(x_{1}, x_{2})}{\lambda}, \forall j = 1, \ldots, deg(g)$. 
Take $\alpha > 0$ such that $\lambda  - \alpha \sigma > 1.$
Thus, taking $(x_{1}, y_{1}) , (x_{2}, y_{2}) \in \T^{d}\times \Sc^{1}$ with $d\big( (x_{1} , y_{1}), (x_{2}, y_{2})\big) < \epsilon_{0},$ for $\epsilon_{0} >0$ small enough,  we have:

$$d\big(F(x_{1}, y_{1}) , F(x_{2}, y_{2})\big) = d_{\T^{d}}(g(x_{1}) , g(x_{2})) + \alpha |f_{x_{1}}(y_{1}) - f_{x_{2}}(y_{2}| = $$
$$\lambda \cdot d_{\T^{d}}(x_{1}, x_{2}) + \alpha \cdot |f_{x_{1}}(y_{1}) - f_{x_1}(y_{2})+f_{x_{1}}(y_2) - f_{x_{2}}(y_{2})| \geq $$
$$\lambda \cdot d_{\T^{d}}(x_{1}, x_{2}) + \alpha \cdot |f_{x_{1}}(y_{1}) - f_{x_{1}}(y_{2})| - \alpha \cdot |f_{x_{1}}(y_2)-f_{x_{2}}(y_2)| \geq $$
$$\lambda \cdot d_{\T^{d}}(x_{1}, x_{2}) + \alpha \cdot |y_{1} - y_{2}| - \alpha \cdot \sigma\cdot|x_{1}-x_{2}| \geq $$
 $$\geq (\lambda - \alpha\cdot \sigma) \cdot d_{\T^{d}}(x_{1}, x_{2}) + \alpha\cdot|y_{1} - y_{2}|.$$
We conclude that $\bf{M}$ is weakly contracting.
\end{proof}

\begin{definition}
     When $\textbf{M}$ is weakly contracting, we define a \emph{natural coupling} $\textbf{P}$ as follows. For each $(x,y)\in \Omega \times \Omega$, we fix permutations $\eta, \sigma$ realizing item (i) of the previous definition. Then for each pair
$(x,y) \in \Omega \times \Omega, t \in \N$ and each word $w=(j_{1},\ldots, j_{t} ) \in \{1,...,k\}^{t}$ we let $\bar{x}^{w}_{t}, \bar{y}^{w}_{t} \in \Omega^{t}$
 be the sequences
$(x_{1},\ldots,x_{t} ),(y_{1},\ldots,y_{t} )$ such that $x_{1} =x^{\eta(j_{1})} \in T^{-1}(x)$ and $x_{1} =x^{\eta(j_{1})} \in T^{-1}(x)$, and for
all $n< t$, $x_{n+1} =(x_{n})^{\eta_{n}(j_{n} )} \in T^{-1}(x_{n})$ and $y_{n+1} =(y_{n})^{\sigma_{n}(j_{n} )} \in T^{-1}(y_{n})$ where $\eta_{n}$ and $\sigma_{n}$ are the permutations associated to the pair $(x_{n},y_{n})$. Then, the natural coupling is
$$
\Pi^{t}_{x,y} = \sum_{w \in \{1, \ldots, k\}}\frac{1}{k^{t}}\delta_{(\bar{x}^{w}_{t} , \bar{y}^{w}_{t})}.
$$
\end{definition}

In other words, the orbits are paired according to the pairing given in the definition of weakly contracting. This ensures that each pair of preorbits remains close.

\begin{definition}\label{def:flat}
    Let $\textbf{M}$ be the uniform backward random walk and $\textbf{P}$ be the natural coupling. Given $0 < \alpha \leq 1$, we say that a
potential $A \in {C}^{\alpha}(\Omega , \R)$ is $\alpha-$\emph{flat} whenever, for some constant $C >0$,
for all $t \in \N$, all $x,y \in \Omega$, and $\Pi^{t}_{x,y}-$almost all $(\bar{x},
\bar{y})$, the following holds:
$$
|A^{t}(\bar{x}) - A^{t}(
\bar{y})| \leq Cd(x,y)^{\alpha};
$$
where 
\begin{align*}
    A^{t} : \Omega^{t} &\xrightarrow[\text{\;\;\;\;\;\;\;\;\;\;\;\;\;\;\;\;\;\;\;\;\;\;\;\;\;\;\;\;\;\;\;\;\;\;\;\;}]{} \R \\
    \bar{x} &\mapsto A(x_{1}) + \cdots + A(x_{t})
\end{align*}
\end{definition}

Then from \cite[Lemma 5.5]{Kl20} comes a flatness criterion:

\begin{lemma}\label{lem:flat}
    Let $\textbf{M}$ be a weakly contracting $1-$to$-k$ transition kernel with contraction
function $c$ and ratio $\lambda >1$. As above, $x^{1}
,\ldots,x^{k}$ denote the points supporting $m_{x}$ and
$\eta, \sigma$ define the natural coupling $\textbf{P}$.
Assume that there is a set $N \subset \Omega$ such that, for all $x,y \in \Omega$ and all $j$, if either $x^{\eta(j)} \notin N$
or $y^{\sigma(j)} \notin N$, then $\ds d(x^{\eta(j)}
,y^{\sigma(j)}) \leq \frac{d(x,y)}{\lambda}$.% In other words, in the definition of "weakly contracting", the indices for which the stronger contraction does not hold are detected by the condition $(x^{\eta(j)},y^{\sigma(j)}) \in N^{2}$.

Let $\alpha \in (0,1]$ and $A \in {C}^{\alpha}(\Omega , \R)$. If for some $C$, all $t$, all $x,y \in \Omega$, and $\Pi^{t}_{
x,y}-$almost all
$(\bar{x},
\bar{y})$ staying in $N$ (i.e. $(x_{n},y_{n}) \in N^{2}$ for all $n \in \{1,\ldots,t\}$) we have
$$
\Big|\sum_{n=1}^{t}A(x_{n})-A(y_{n})\Big| \leq Cd(x,y)^{\alpha}, 
$$
then $A$ is $\alpha-$flat.
\end{lemma}

Finally, \cite[Theorem 5.8]{Kl20} provides a result that guarantees the spectral gap property.

\begin{theorem}\label{theor:weagap}
    Suppose that the uniform backward random walk $\textbf{M}$ is weakly contracting and that the potential $A \in {C}^{\alpha}(\Omega , \R) $ is $\alpha-$flat. Then the transfer operator $\mathcal{L}_{F, \phi}$ has the spectral gap property, acting on ${C}^{\alpha}(\Omega , \C)$.
\end{theorem}

\

\begin{proof}[Proof of Theorem \ref{maintheoC}]
Define $\tilde{\mathcal{H}}:= \{\phi \in C^{r}(\mathbb{T}^{d} \times \Sc^{1} , \R) : \exists \epsilon > 0 \text{ with } \phi(x , y) = \phi(x , q_{j,x}) %\text{ and } \phi(\tilde{x}, \tilde{y}) = \phi(x_{i}, \tilde{y}) 
\text{ for all } x
%, \tilde{y}
\in \mathbb{T}^{d}, |y - q_{j,x}| < \epsilon
%d(\tilde{x}, x_{i}) < \epsilon, i = 1, \ldots, l-1 
\text{ and } j = 1, \ldots, k \}
$. Note that $\tilde{\mathcal{H}}$ is dense in $C(\mathbb{T}^{d} \times \Sc^{1} , \R)$ and $t \tilde{\mathcal{H}} \subset \tilde{\mathcal{H}},$ for all $t \in \R.$ We will show that if $\phi \in \tilde{\mathcal{H}}$ and $0 < \alpha < 1 $ then $\mathcal{L}_{F , \phi}$ has the spectral gap property, acting on $C^{\alpha}(\mathbb{T}^{d} \times \Sc^{1} , \C);$ it's enough to prove Theorem \ref{maintheoC}.

Take $\tilde{F} := F^{m},$ for $m$ big enough.  Note that there exists $0 < \sigma < 1$ such that: given $z_{1}, z_{2} \in \mathbb{T}^{d} \times \Sc^{1}$, with $z_{1} \neq z_{2},$ we have $\tilde{F}^{-1}(z_{1}) = \{w_{1,1} , \ldots, w_{1, deg(F)^{m}}\}$ and $\tilde{F}^{-1}(z_{2}) = \{w_{2,1} , \ldots, w_{2, deg(F)^{m}}\}$ with $d(w_{1,1}, w_{2,1}) \leq \sigma \cdot d(z_{1}, z_{2})$ and $d(w_{1,j}, w_{2,j}) < d(z_{1}, z_{2}),$ for all $j = 2, \ldots, deg(F)^{m}.$ Taking the uniform backward random walk $\textbf{M} = (m_{z})_{z \in \mathbb{T}^{d} \times \Sc^{1}},$ where $\ds m_{z} := \frac{1}{deg(F)^{m}}\sum_{F^{m}(w) = z}\delta_{w},$ we have that $\textbf{M}$ is an  $1-$to$-deg(F)^{m}$ weakly contracting  transition kernel (according to Definition \ref{def:transker}). 
Note that fixed $\epsilon > 0$ there exists  $\lambda > 1$ such that $||D\tilde{F}(x , y)|| \geq \lambda$, if 
$\ds\inf_{j = 1, \ldots, k }|y - q_{j,x}| \geq \frac{\epsilon}{2}$.% and $\ds\inf_{i = 1, \ldots, l-1}d(x, x_{j}) \geq \frac{\epsilon}{2}$.

Thus, given $\phi \in \tilde{\mathcal{H}}$ and $d(z_{1}, z_{2}) < \frac{\epsilon}{2}$ 
take $$N := \Big\{(x , y ) \in \mathbb{T}^{d} \times \Sc^{1} :  |y - q_{j,x}| < \frac{\epsilon}{2}, %d(\tilde{x}, x_{i}) < \frac{\epsilon}{2},
\text{ for some }  j = 1, \ldots, k %\text{ and } i = 1, \ldots, l-1 
\Big\},
$$
according to Lemma \ref{lem:flat}.

Note that $\tilde{F}^{-1}(z_{1}) = \{w_{1,1} , \ldots, w_{1, deg(F)^{m}}\}$ and $\tilde{F}^{-1}(z_{2}) = \{w_{2,1} , \ldots, w_{2, deg(F)^{m}}\}$ such that for each $w_{1 ,l} = (x_{1, l} , y_{1,l})$ and $w_{2,l} = (x_{2,l}, y_{2,l})$ with $w_{1 ,l} , w_{2 ,l} \in N$ we have that $|y_{1,l} - q_{j,x_{1,l}}| < \frac{\epsilon}{2} , |y_{2,l} - q_{j,x_{2,l}}| < \frac{\epsilon}{2},$ for some $j \in \{1,\ldots, k\}$. In this case, $$|\phi(w_{1,l}) - \phi(w_{2,l})| = |\phi(x_{1,l}, q_{j,x_{1,l}}) - \phi(x_{2,l}, q_{j,x_{2,l}})| \leq 
$$
$$||\phi||_{1} \cdot \Big[ d(x_{1,l} , x_{2,l}) + \alpha \cdot |q_{j,x_{1,l}} - q_{j,x_{2,l}}| \Big]
\leq  ||\phi||_{1} \cdot (1 + C) d(x_{1,l} , x_{2,l}).$$
 
Suppose that $$\tilde{F}^{-n}(z_{1}) = \{w_{n, 1,1} , \ldots, w_{n , 1, deg(F)^{nm}}\}
 \text{ and } \tilde{F}^{-n}(z_{2}) = \{w_{n, 2,1} , \ldots, w_{n , 2, deg(F)^{nm}}\},$$ with $w_{i, 1, \eta(i)}, w_{i, 2, \eta(i)} \in N$, for all $i = 0 , \ldots, n.$ Then 
$$\Big|\sum_{i =1}^{n}\phi(w_{i, 1, \eta(i)}) - \phi(w_{i, 2, \eta(i)})\Big| 
\leq ||\phi||_{1} \cdot (1+ C)\sum_{i =1}^{n}d(x_{1,i} , x_{2,i}) \leq 
$$
$$||\phi||_{1} \cdot (1+ C)\sum_{i =1}^{n}\frac{1}{\inf_{x \in \mathbb{T}^{d}} ||Dg(x)||^{i}}d(x_{1,0},x_{2,0})\leq 
(1+ C) \frac{\inf_{x \in \mathbb{T}^{d}} ||Dg(x)||}{\inf_{x \in \mathbb{T}^{d}} ||Dg(x)|| - 1}\cdot d(z_{1}, z_{2}).$$

Applying the Lemma \ref{lem:flat}, we have that $\phi$ is $\alpha-$flat (according to Definition \ref{def:flat}), for all $0 < \alpha < 1$. It follows from the 
Theorem \ref{theor:weagap}
that $\mathcal{L}_{\tilde{F}, \phi}$ has the spectral gap property, acting on $C^{\alpha}(\mathbb{T}^{d} \times \Sc^{1} , \C).$

Finally, we will show that $\mathcal{L}_{F , \phi|C^{\alpha}}$ has the spectral gap property for all $\phi \in \tilde{\mathcal{H}}.$ Note that $\mathcal{L}_{F , \phi}^{m} = \mathcal{L}_{\tilde{F}, S_{m}\phi},$ where $S_{m}\phi := \sum_{j=0}^{m-1}\phi \circ F^{j}$.
Since $\phi \in \tilde{\mathcal{H}}$ then $S_{m}\phi \in \tilde{\mathcal{H}}$, thus $\mathcal{L}_{F , \phi|C^{\alpha}}^{m}$ has the spectral gap property. Since $\mathcal{L}_{F, \phi}$ is a positive operator then there exists a probability $\nu$ such that $\mathcal{L}_{F , \phi}^{\ast}\nu = \rho(\mathcal{L}_{F , \phi|C^{0}})\nu$. On the other hand, the spectral gap property of $\mathcal{L}_{F , \phi|C^{\alpha}}^{m}$ implies that $\rho(\mathcal{L}_{F , \phi|C^{0}})^{m} = \rho(\mathcal{L}_{F , \phi|C^{\alpha}}^{m})$ and $\ds\lim_{n \to +\infty} \mathcal{L}_{F , \phi|}^{n}g = \lim_{n \to +\infty} \mathcal{L}_{F , \phi|}^{mn}g = \int g d\nu \cdot h_{\phi}$, for all $g \in C^{\alpha}(\mathbb{T}^{d} \times \Sc^{1} , \C)$. We conclude then that $\mathcal{L}_{F , \phi|C^{\alpha}}$ has the spectral gap property. 
\end{proof}

\section{Questions and comments}\label{sec:quest}

In our Theorem \ref{theorA}, we assume that the expansive dynamic $f : M \rightarrow M$ admits a repeller periodic point; thus, a natural question is about the necessity of the hypotheses.

\begin{mainquestion}
    If a $C^{1}-$local diffeomorphism $f : M \rightarrow M$ is expansive, then $f$ has a repeller periodic point ?
\end{mainquestion}

Another natural question concerns additional information in Theorem \ref{theorA} about the spectral gap property and the non-existence of phase transitions in all parameters.

\begin{mainquestion}
    Is the subset $\Big\{\phi \in C^{r}(M , \R) : \phi \text{ has no phase transition }\Big\}$ dense ?
\end{mainquestion}

\begin{mainquestion}
    Is the subset 
    $$\Big\{\phi \in C^{r}(M , \R) : \mathcal{L}_{t\phi|C^{r}} \text{ has  the spectral gap property for all } t \in \R\Big\}$$ dense ?
\end{mainquestion}

%%%%%%%%%%%%%%%%%%%%%%%%%

\vspace{.3cm}
\subsection*{Acknowledgements.}
 TB was partially supported by “Conselho Nacional de Desenvolvimento Científico e Tecnológico – CNPq” (Grants PQ-2021 and Universal-18/2021).
%The authors are deeply grateful to ... for valuable comments.

%%%%%%%%%%%%%%%%%%%%%%%%%%%%%%%%%%%%
\bibliographystyle{alpha}

\end{document}